\documentstyle[12pt]{article}

\font\teneufm=eufm10
\font\seveneufm=eufm7
\font\fiveeufm=eufm5
\newfam\eufmfam
\textfont\eufmfam=\teneufm
\scriptfont\eufmfam=\seveneufm
\scriptscriptfont\eufmfam=\fiveeufm

\newfam\msbfam
\font\tenmsb=msbm10 scaled \magstep1  \textfont\msbfam=\tenmsb
\font\sevenmsb=msbm7 scaled \magstep1 \scriptfont\msbfam=\sevenmsb
\font\fivemsb=msbm5 scaled \magstep1  \scriptscriptfont\msbfam=\fivemsb
\def\Bbb{\fam\msbfam \tenmsb}

\def\RR{{\Bbb R}}
\def\CC{{\Bbb C}}

\def\PP{{\Bbb P}}
\def\NN{{\Bbb N}}
\def\ZZ{{\Bbb Z}}

\def\Aut{\hbox{Aut}}
\def\Re{\hbox{Re}\,}
\def\Im{\hbox{Im}\,}
\def\ra{\rightarrow}
\def\zhat{\widehat{z}}
\def\zbar{\overline{z}}

\def\e{\epsilon}

 \def\HollowBoxx #1#2#3{{\dimen0=#1 \advance\dimen0 by -#2       
       \dimen1=#1 \advance\dimen1 by #3                       
        \vrule height 0pt depth #3 width #2                   
       \hskip -#3
       \vrule height #1 depth #3 width #3}}                   
 \def\LeftContraction{\mathord{\kern1.45pt \HollowBoxx{6pt}{3.5pt}{.4pt}}\,}

 \def\HollowBox #1#2#3{{\dimen0=#1 \advance\dimen0 by -#3       
       \dimen1=#1 \advance\dimen1 by #3                       
        \vrule height #1 depth #3 width #3                    
        \vrule height 0pt depth #3 width #2                   
        \hskip -#3}}                                             
 \def\RightContraction{\mathord{\, \HollowBox{6pt}{3.1pt}{.4pt}} \kern1.6pt}

\newtheorem{theorem}{THEOREM}[section]

\newtheorem{corollary}[theorem]{Corollary}

\newtheorem{example}[theorem]{EXAMPLE}

\begin{document}
\begin{center}
{\Large \bf Domains with Non-Compact} 
\medskip \\
{\Large \bf Automorphism Group:  A Survey}  
\medskip \\
{\normalsize \rm A. V. Isaev \ \ \ and \ \ \ S. G. Krantz}
\end{center} 

\begin{quotation} 
\small \sl
We survey results arising from the study of domains in
$\CC^n$ with non-compact automorphism group.  Beginning
with a well-known characterization of the unit ball,
we develop ideas toward a consideration of weakly pseudoconvex 
(and even non-pseudoconvex) domains with particular emphasis on 
characterizations
of {\bf (i)} smoothly bounded domains with non-compact
automorphism group and {\bf (ii)}  the Levi geometry of
boundary orbit accumulation points.

Particular attention will be paid to results derived in the 
past ten years.
\end{quotation}

\markboth{A.\ V.\ Isaev and S.\ G.\ Krantz}{Domains with 
   Non-compact Automorphism Group}

\setcounter{section}{-1}

\section{Introduction}

In any area of mathematics, one of the fundamental problems is 
to determine the equivalence of the structures under 
consideration---that is, to determine the morphisms in
the relevant category.  In complex
analysis one is interested, for example, in the holomorphic
equivalence of complex manifolds. The problem that we
study here is somewhat more subtle:  we wish to see to
what extent a domain is determined by the group of its biholomorphic
self-mappings.

In this paper we deal only with domains in ${\CC}^n$, for even then
the equivalence problem that we wish to study (described below) turns
out to be highly non-trivial. It is a well-known fact that, if $n\ge
2$, then two given domains in $\CC^n$ are most likely to be
holomorphically inequivalent. This can be understood, for example, by
examining the induced mapping between the boundaries of two given
domains (in cases when the original mapping can be extended to a
mapping between the closures of the domains). Poincar\'e was one of
the first to notice the connection between the equivalence problem
for domains and that for their boundaries. Considering domains with
real analytic boundary in $\CC^2$ and writing the equations of the
boundaries in a special form, he showed that Taylor expansions that
have different coefficients for monomials of sufficiently high degree
define inequivalent boundaries. Thus
Poincar\'e proved that there
are infinitely many inequivalent domains \cite{Po}.

As shown in \cite{Fe}, any biholomorphic mapping between smoothly
bounded strictly pseudoconvex domains {\it does} extend to a mapping
between their closures, and therefore one can endeavor to find an
analogue of Poincar\'e's argument in this case. It is possible, for
example, to derive such an analogue for general real analytic 
strictly pseudoconvex hypersurfaces (as well as for any
hypersurfaces with non-degenerate Levi form) from Moser's normal form
for their defining functions \cite{CM}. Moreover, it turns out that
almost any two strictly pseudoconvex domains with only $C^2$-smooth
boundary are inequivalent \cite{GK2} (see also
\cite{GK1}, \cite{BSW}).

We have gone into some considerable detail on this point
of generic domain inequivalence in order to emphasize the
special nature of the function theory of several complex
variables, and to stress the particular difficulties that
we shall face.  Note especially that there is no moduli space,
nor anything like a Teichm\"{u}ller space, for smoothly
bounded domains in $\CC^n$ (in fact this assertion has
been proved in a strong sense in \cite{LR}).

As a result of these considerations, we must restrict ourselves to special
collections of domains that on the one hand are sufficiently small
so that we may hope for a reasonable classification, and on the
other hand are sufficiently large to possess a rich and interesting
structure.

Let $D\subset {\CC}^n$ be a bounded (or, more generally, 
Kobayashi-hyperbolic---see Section 1 for the definition) domain. Denote by
$\hbox{Aut}(D)$ the group (under composition) of holomorphic
automorphisms of $D$. The group $\hbox{Aut}(D)$ is a topological
group with the natural topology of uniform convergence on compact
subsets of $D$ (the compact-open topology). It turns out that
$\hbox{Aut}(D)$ can be given the structure of a Lie group whose
topology agrees with the compact-open topology (see \cite{Kob1}).
Many abstract Lie groups can be realized as the automorphism groups
of bounded domains in complex space \cite{SZ}, \cite{BD}, \cite{TS},
but in this paper we deal only with domains for which $\hbox{Aut}(D)$
is ``large enough''. 

More precisely, we consider the class of domains for which
$\hbox{Aut}(D)$ is {\it non-compact}. By a classical theorem of H. Cartan
(see \cite{N}), for a bounded domain this condition is equivalent to the
non-compactness of every orbit of the action of $\hbox{Aut}(D)$ on $D$
(which is in fact equivalent to the existence of only one non-compact
orbit). For example, any homogeneous domain (i.e. domain on which
$\hbox{Aut}(D)$ acts transitively) has non-compact automorphism group.
The study of bounded homogeneous domains was initiated by \'E. Cartan
\cite{Car} and eventually led to their complete classification
\cite{P-S} (for the more general case of complex spaces see e.g.
\cite{HO}). The technique by which these classifications were
obtained is mostly algebraic.  We will, however, be more interested
in geometric and analytic methods that have been developed under the
additional hypothesis of a certain regularity of the boundary of the
domain (local or global $C^{\infty}$-smoothness in many cases). The
regularity of the boundary is indeed a crucial component of all the
considerations below; to illustrate its importance, we only mention
here that for homogeneous domains with $C^2$-smooth boundary the
classification in \cite{P-S} turns into a single domain, namely the
unit ball.  We also note that, when the boundary smoothness is less
than $C^2$, then many of the basic ideas in this subject break down
(see \cite{GK2}).

Section 1 contains basic definitions, notation and elementary
background material. The reader already familiar with this subject may
safely skip Section 1 and refer back to it as needed. We begin our survey 
in Section 2 with the now classical Ball Characterization
Theorem for strictly pseudoconvex domains; this is the first main
result in the subject (from the point of view that we wish to
promulgate); it, in turn, led to the Greene/Krantz theorems (and
their generalizations) that were the first attempts to obtain a
general result for weakly pseudoconvex domains. Section 3
is built around the domains of Bedford/Pinchuk. It is quite
a plausible conjecture that these domains give a complete classification
of smoothly bounded domains with non-compact automorphism group.  

The techniques in Section 3 clearly show the importance of the
hypothesis of finiteness of type in the sense of D'Angelo
\cite{D'A1} of the boundary of the domain at the boundary orbit
accumulation points (by definition, a point
$q\in\partial D$ is a boundary orbit accumulation point for the
action of $\hbox{Aut}(D)$ on $D$ if there exist a point $p\in D$ and
a sequence $\{f_j\}\subset \hbox{Aut}(D)$ such that
$f_j(p)\rightarrow q$ as $j\rightarrow\infty$). The conjecture that
finiteness of type always obtains at the boundary orbit accumulation points of
a smoothly bounded domain is known as the ``Greene/Krantz
conjecture'' and is discussed in Section 4. Another hypothesis
important for many results in Section 3 is the pseudoconvexity of the boundary
near a boundary orbit accumulation point. It is also discussed in
Section 4. 

In Section 5 we deal with properties of the boundary orbit
accumulation set (the set of all boundary orbit accumulation points)
as a whole,
in particular, certain extremal properties of some invariants of the
boundary of the domain.  In Section 6 we consider domains with less
than $C^{\infty}$-regularity of the boundary (e.g. finitely smooth or
piecewise smooth) and also unbounded domains. Some of the results for
unbounded domains are localizations of those mentioned in the
preceding sections, but some of them are essentially global, and the
domain is then required to be Kobayashi-hyperbolic. For domains with rough
boundary, the results included in this section lead, in particular,
to an analogue of the Bedford/Pinchuk domains in the finitely smooth
case.  Note again that, when the domain under study has {\it extremely
rough boundary}---say fractal boundary---then the classification
problem appears to be intractable (see \cite{Kra4}).

To set the tone for this article, we now present four examples of
domains in $\CC^2$ which, taken together, tend to suggest some of the
subtlety and beauty of the subject.  They are the domains
\begin{eqnarray*}
B^2 &:=& \{(z_1,z_2): |z_1|^2 + |z_2|^2 < 1\}, \\
E_{1,2} &:=& \{(z_1,z_2): |z_1|^2 + |z_2|^4 < 1\}, \\
E_{2,2} &:=& \{(z_1,z_2): |z_1|^4 + |z_2|^4 < 1\}, \\
E_{1,\infty} &:=& \{(z_1,z_2): |z_1|^2 + 2
\cdot e^{-1/|z_2|^2} < 1\}.
\end{eqnarray*}

The domain $B^2$ is the unit ball, and has transitive automorphism
group---this follows, for example, from the explicit description of
$\hbox{Aut}(B^2)$ (see \cite{Ru} or \cite{Kra3}).

The domain $E_{1,2}$ has non-compact automorphism group; to wit,
the automorphisms
$$
(z_1,z_2) \longmapsto \left ( \frac{z_1 - a}{1 - \overline{a}z_1},
              \frac{\sqrt[4]{1 - |a|^2}z_2}{\sqrt{1 - \overline{a}z_1}} \right ) \ ,
\qquad |a| <1 \ ,
$$
form a non-compact set of automorphisms of $E_{1,2}$. As the parameter
$a$ approaches $-1$, the above family moves any interior point of
$E_{1,2}$ to the boundary point $(1,0)$ which is therefore a boundary
orbit accumulation point for $\hbox{Aut}(E_{1,2})$. The automorphism
group of $E_{1,2}$ cannot be transitive, because then $E_{1,2}$ would
be holomorphically equivalent to the unit ball (by the Ball
Characterization Theorem--Theorem 2.1
below), but it is not (for instance, by a theorem of Bell \cite{Bel1}).

The domain $E_{2,2}$ has {\it compact}\/ automorphism group.  The assertion
is not entirely obvious.  This example differs from the preceding one
in that the boundary of $E_{2,2}$ has {\it two} orthogonal circles
of weakly pseudoconvex points (see Section 1 for the definition),
while $E_{1,2}$ has just one.  Any automorphism of $E_{2,2}$ must
{\bf (i)}  extend smoothly to the boundary (see e.g. \cite{Bel1})
and {\bf (ii)}  take weakly pseudoconvex points to weakly pseudoconvex
points.  The two circles of weakly pseudoconvex points therefore
serve to harness any fixed compact subset of $E_{2,2}$; in particular,
they prevent any orbit from accumulating at a point in the boundary.
Thus, the automorphism group of $E_{2,2}$ is compact. 

The domain $E_{1,\infty}$ also has compact automorphism group.  This is the
most subtle example of all.  The questions that it raises will be the
focus of much of the rest of the present survey.  Briefly, $E_{1,\infty}$
has compact automorphism group for the following reason. If the automorphism
group is non-compact then some orbit must accumulate at the boundary.
If it accumulates at a point of the form $q=(q_1,q_2)$ with $q_2 \ne 0$
then $q$ is a point of strong pseudoconvexity (see Section 1 for the definition).  It then
follows from the Ball Characterization Theorem that
$E_{1,\infty}$ is holomorphically equivalent to the unit ball---which it is not.
If instead the orbit accumulates at a point of the form
$q=(q_1,0)$, then $q$ is infinitely
flat in a sense to be made precise in the next section.  It turns
out (and there is a general conjecture to this effect---see Section 4) that
infinitely flat points in this sense {\it cannot} be boundary orbit accumulation
points. Further details on this example can be found in \cite{GK7}.
An alternative proof of compactness of $\hbox{Aut}(E_{1,\infty})$
follows from Theorem 3.6 in Section 3 where we discuss the case
of {\it Reinhardt domains}. Note that the domains
$E_{2,2}$ and $E_{1,\infty}$ are {\it not} holomorphically equivalent
(this follows, for instance, from \cite{BKU}). 

We have exhibited four domains of the same topological
type---indeed the closure of each one is diffeomorphic to the closure
of each of the others---yet with
strikingly different holomorphic automorphism group
characteristics (note that
these domains are pairwise
holomorphically inequivalent). The reasoning that we have sketched in this brief
discussion sets the tone for the arguments that we shall present
in the rest of this paper.

Before proceeding, we note that the reader may find it useful to
compare the present paper with the earlier surveys \cite{GK5},
\cite{Kra2}, \cite{Kra4}.

We would like to thank E.\ Bedford, J.\ D'Angelo, H.\ Gaussier, R.\ Greene
and K.-T.\ Kim for their interest in our work and for useful suggestions.

\section{Preliminaries}

In this section we give the definitions of the main concepts and some
of the facts that we use later in the paper. We shall not discuss
them here in any detail, nor shall we make any historical remarks; we
refer the reader to \cite{Kra3} for additional information and
background.    

Although it may be possible to profitably study manifolds
with non-compact automorphism group, in
this paper we restrict attention to domains in $\CC^n$,
where by a {\it domain} we mean a connected open set.  Our domains
are usually, but not always, bounded; we generally denote them
by $D$ or $\Omega$.

A {\it holomorphic function} on a domain $D$ is a function
\begin{eqnarray*}
f: D & \rightarrow & \CC \\
    (z_1,\dots,z_n) & \mapsto & f(z_1,\dots,z_n)
\end{eqnarray*}
that is holomorphic in each variable separately.  Such a function
is automatically $C^\infty$-smooth as a function of the $2n$ real
variables $x_1, y_1, \dots, x_n, y_n$, where
$z_j = x_j + i y_j$, $j = 1,\dots, n$.  Such a function also
has a locally convergent power series expansion near every point of $D$.
If $D_1$, $D_2$ are domains in $\CC^n$ then a {\it holomorphic
mapping} of $D_1$ to $D_2$ is a mapping
$$
F:D_1 \rightarrow D_2
$$
such that $F(z_1,\dots,z_n) = \bigl ( f_1(z_1,\dots,z_n), \dots
f_n(z_1,\dots,z_n) \bigr )$ where each $f_j$ is a holomorphic function.
We say that $F$ is {\it biholomorphic} if it is one-to-one and
onto.  A biholomorphic mapping has an inverse which is automatically
biholomorphic itself. Two domains are called {\it holomorphically
equivalent} if there is a biholomorphic mapping from one domain onto
the other.

The collection of biholomorphic mappings of a domain $D$ onto itself
(often termed {\it biholomorphic self-maps} or {\it automorphisms}) 
of $D$ clearly forms a group under composition of mappings.  We
denote this group by $\hbox{Aut}(D)$.  This group is given the
topology of uniform convergence on compact subsets of $D$ (the
compact-open topology). So equipped, $\hbox{Aut}(D)$ turns out to be
a real Lie group when $D$ is bounded.  If the group is positive dimensional,
then it is never a complex Lie group for any bounded $D$ (see
\cite{Kob1}). 

In the complex plane, a bounded domain with $C^1$-smooth boundary 
that is finitely connected (with
connectivity at least one) has compact automorphism group---in
fact if the connectivity is at least two then the automorphism
group is finite (see \cite{GK5} for details).  The only
bounded planar domain with $C^1$-smooth boundary and non-compact automorphism group is, up to a
biholomorphism, the unit
disc (see \cite{Kra1}).

As we will see below, in dimensions 2 and higher, the collection of
domains with non-compact automorphism group and regular boundary is much bigger.
For the most part, our discussion will center on domains with
{\it smooth boundary}. Let $D \subset \CC^n$ be a domain. In a
neighborhood $U$ of a fixed point $p\in\partial D$ we can
write 
$$ 
D\cap U = \{z \in U: \rho(z) < 0\}. 
$$
Such a
function $\rho$ is called a {\it defining function for $D$ near
$p$}.  We say that, for $1\le k\le\infty$, $D$ has $C^k$-smooth or
real analytic boundary near $p$ if there is a defining function
$\rho$ for $D$ near $p$ which is, respectively, either $C^k$-smooth 
or real analytic and  $\nabla \rho \ne 0$ on
$\partial D$. The boundary is said to be globally $C^k$-smooth or 
real analytic if it is such at every point.
When the boundary is globally $C^k$-smooth then it is easy to patch
together local defining functions to obtain a single global defining
function for the entire boundary. From now on, when speaking about defining functions of domains with
smooth boundary, we will be assuming that these functions satisfy the
conditions just discussed.

It is natural in our studies to pay special attention to {\it
pseudoconvex domains} or {\it domains of holomorphy}. A domain $D$
is called a domain of holomorphy if it is the natural
domain of existence for some holomorphic function; in other words,
if there exists a function holomorphic in $D$ and such that it cannot
be holomorphically continued past any boundary point of $D$. A more
technical equivalent definition is as follows: a domain $D$
is a domain of holomorphy if, 
for any compact set
$K\subset D$, its {\it holomorphic hull} $\hat K:=\{z\in
D:|f(z)|\le\hbox{max}_{\zeta\in K}|f(\zeta)|,\,\hbox{for any $f$
holomorphic in $D$}\}$ is also compact in $D$. 

If $D$ has at least
$C^2$-smooth boundary near $p\in\partial D$, then $\partial D$ is said
to be {\it pseudoconvex at} $p$ if there is a defining
function $\rho$ for $D$ near $p$ such that
$$
{\cal L}_{\rho}(p)(w,w):=\sum_{j,k = 1}^n \frac{\partial^2 \rho}
{\partial z_j \partial \overline{z}_k}(p)
w_j \overline{w}_k \ge 0 \eqno{(1.1)}
$$
for all $w:=(w_1,\dots,w_n)\in T^c_p(\partial D)$;
here $T^c_p(\partial D)$ is the {\it complex 
tangent space to $\partial D$ 
at $p$}, which is the
maximal complex subspace of the ordinary real tangent space 
$T_p(\partial D)$:
$$
T^c_p(\partial D):= \left \{(w_1,\dots,w_n):
\sum_{j=1}^n \frac{\partial \rho}{\partial z_j}(p) w_j = 0 \right \} .
$$
We call $p\in\partial D$ a point of {\it strong} or {\it strict}
pseudoconvexity if the inequality (1.1) is strict for non-zero $w\in
T^c_p(\partial D)$. The Hermitian form ${\cal L}_{\rho}(p)$ defined
in (1.1) is called the {\it Levi form of $\partial D$ at $p$ }. It
depends on the defining function $\rho$ and is defined up to
multiplication by a positive constant; therefore, the signs of the
eigenvalues of ${\cal L}$ do {\it not} depend on the choice of
$\rho$.  In particular, the notions of pseudoconvexity and strict
pseudoconvexity at $p$ do {\it not} depend on $\rho$. 

It turns out that a domain with $C^2$-smooth boundary is a domain of
holomorphy if and only if its boundary is pseudoconvex at every
point.
A pseudoconvex domain whose boundary is strictly pseudoconvex
at every point is called {\it strictly pseudoconvex}.

The other extremal situation---in contrast with strict
pseudoconvexity---is when the Levi form is identically zero in a
neighborhood of $p\in\partial D$. In this case $\partial D$ is called
{\it Levi-flat near $p$} and is then foliated near $p$ by complex
submanifolds of dimension $n-1$; conversely, if $\partial D$ admits
such a foliation, it is Levi-flat (see e.g. \cite{T}).

We next turn to the notion of {\it type} in the sense of D'Angelo 
for
$C^{\infty}$-smooth real hypersurfaces in $\CC^n$ \cite{D'A1}. 
The
type measures the order of tangency of (possibly singular)
holomorphic curves with the hypersurface at a given point. Let
$D\subset \CC^n$ be a domain with $C^{\infty}$-smooth boundary and
let $p\in \partial D$. Then the type $\tau(p)$ of $\partial D$ at $p$ is
defined as
$$
\tau(p):=\sup_{F}\frac{\nu(\rho\circ F)}{\nu(F)},
$$
where $\rho$ is a defining function of $D$ near $p$, the supremum is
taken over all holomorphic mappings $F$ defined in a neighborhood of
$0\in\CC$ into $\CC^n$ such that $F(0)=p$, and $\nu(\phi)$ is the
order of vanishing of a function $\phi$ at the origin. The boundary
$\partial D$ is said to be of {\it finite type at $p$} if
$\tau(p)<\infty$. The domain $D$ is a {\it domain of finite type} if
$\partial D$ is of finite type at every point. It is an important
fact that, if $D$ is a bounded domain of finite type, then the type is
uniformly bounded on $\partial D$; this last fact follows from a weak
semi-continuity property of $\tau$ (see \cite{D'A2}).  
Examples of domains of
finite type are bounded domains with real analytic boundary
\cite{D'A2}, \cite{DF}, \cite{L}---though in many ways these
examples are atypically simple. We also note that the boundary of a
domain is of finite type at the points of strict
pseudoconvexity. Occasionally we will be using a weaker
condition than that of finite type: we say that $\partial D$ is
{\it variety-free at} $p\in\partial D$
if $\partial D$ does not
contain positive-dimensional complex varieties passing through $p$.

Domains of finite type are important for function theory.  They have
many of the attractive properties of strongly pseudoconvex domains; in
particular, biholomorphic mappings
of bounded domains of finite type extend to diffeomorphisms of the closures.

In this paper we mainly consider bounded domains. However, we will see
that some of the results and techniques can be generalized to {\it
Kobayashi-hyperbolic} domains. Hyperbolicity is geometrically a
natural generalization of boundedness and is defined in terms of the
{\it Kobayashi pseudometric}. Let $M$ be a complex manifold, $p\in M$,
$v\in T_p(M)$. The {\it Kobayashi pseudonorm} of $v$ is
the quantity
$$
k(p,v):=\inf_{F}\left\{\frac{1}{r}\right\},
$$
where the infimum is taken over all holomorphic mappings $F$ from discs
$\Delta_r:=\{z\in\CC: |z|<r\}$ to $M$ such that $F(0)=p$, $F'(0)=v$. For a
connected $M$ the Kobayashi pseudometric $K(p,q)$, $p,q\in M$, can now be
defined as
$$
K(p,q):=\inf_{\gamma}\int_{0}^{1}k(\gamma(t),\gamma'(t))\, dt ,
$$
where the infimum is taken over all smooth paths
$\gamma:[0,1]\rightarrow M$ that join $p$ and $q$ \cite{PoSh}. The
Kobayashi pseudometric is a biholomorphic invariant and generalizes the
Poincar\'e metric on the unit disc in $\CC$. 

A manifold $M$ is called
Kobayashi-hyperbolic if the Kobayashi pseudometric on $M$ is in fact a
metric. To illustrate that hyperbolicity is really a generalization
of boundedness we mention here that hyperbolic manifolds possess the
Liouville property which clearly holds for bounded domains: a
holomorphic mapping from $\CC$ into a hyperbolic manifold must be
constant (see \cite{Kob1} for an elegant discussion of the relation
of hyperbolicity, the Liouville property, and curvature).
A complex manifold $M$ is called {\it complete hyperbolic}, if it is
hyperbolic and, in addition, the Kobayashi metric on $M$ is
complete. Examples of complete hyperbolic manifolds are bounded
strictly pseudoconvex domains in $\CC^n$ \cite{Gr}. 

Another invariant metric (which is going to be a Hermitian metric)
that we mention here is the {\it Bergman metric}. Let $D\subset
\CC^n$ be a domain. Let $\{\phi_j\}_{j=1}^{\infty}$ be an
orthonormal basis in the space of
holomorphic square integrable functions on $D$.  The function
$$
B(p,q):=\sum_{j=1}^{\infty}\phi_j(p)\overline{\phi_j(q)}, \qquad p,q\in
D
$$
is called the {\it Bergman kernel} of $D$. The Bergman metric is then
defined as
$$
ds^2_B:=\sum_{j,k=1}^n\frac{\partial^2 B(z,z)}{\partial
z_j\,\partial{\overline z_k}}dz_j\,d\overline{z_k}.
$$
  
We will also need some invariant volume elements. Let $D\subset\CC^n$
be a domain and $p\in D$. The {\it Carath\'eodory volume element of
$D$ at $p$} is defined to be
$$
V^C(p):=\sup_{F}\{|\det F'(p)|\},
$$
where the supremum is taken over all holomorphic mappings $F$ from $D$ to
the unit ball $B^n:=\{(z_1,\dots,z_n):|z_1|^2+\dots+|z_n|^2<1\}$ such
that $F(p)=0$. Likewise, the {\it Eisenman-Kobayashi volume element of $D$
at $p$} is defined by
$$
V^K(p):=\inf_{F}\left\{\frac{1}{|\det F'(0)|}\right\},
$$
where the infimum is taken over all holomorphic mappings $F$ from $B^n$
into $D$ such that $F(0)=p$. The quotient $\frac{V^C}{V^K}$ is a
biholomorphic invariant and will be called the {\it C/K-invariant}
(see \cite{Kra3}, \cite{GK5} for a detailed discussion of this invariant
and its uses). 

\section{The Ball Characterization Theorem and  \hfill \break
Theorems of Greene/Krantz Type}

The first result that we mention in this section is the famous Ball
Characterization Theorem of Bun Wong and Rosay:

\begin{theorem}[\cite{Ro}]  \sl  Let $D\subset{\CC}^n$ be a bounded
domain with $Aut(D)$ non-compact.  Assume that there exists a
boundary orbit accumulation point in a neighborhood of which
$\partial D$ is
$C^2$-smooth and strictly
pseudoconvex. Then $D$ is holomorphically equivalent to the unit ball
$B^n$.
\end{theorem}

This result was first proved in \cite{W1} for globally strictly
pseudoconvex domains. An alternative proof (in the case of
$C^{\infty}$-smooth boundary) based on an analysis of
the holomorphic sectional curvature of the Bergman metric was obtained
in \cite{Kl} (for related results see also \cite{KY}). 

It is important here to realize that the Levi geometry of the
boundary orbit accumulation point completely determines the entire domain.
Thus (micro)local geometric information at the boundary orbit accumulation
point gives global geometric information.  This theme will be
one of the unifying ideas in the remainder of the present paper.
The original approach of Bun Wong (to construct a special metric
using the hypotheses), the Bergman geometry approach of Klembeck,
and the function-theoretic approach of Rosay all manifest this
local/global dialectic in different ways.

Theorem 2.1 clearly implies the following alternative characterization
of $B^n$ (cf. \cite{P-S}):

\begin{corollary} \sl  If $D\subset{\CC}^n$ is a bounded homogeneous
domain with $C^2$-smooth boundary, then $D$ is biholomorphically
equivalent to $B^n$.  
\end{corollary}
{\bf Proof:}  Any $C^2$-smooth
bounded domain in Euclidean space has a
boundary point that is strongly convex, hence strongly pseudoconvex
(just take a fixed point $p$ in space that is far away from the
domain and then take a point in the boundary of $D$ that is at the
maximal distance from $p$).  Now apply Theorem 2.1.  \hfill $\Box$
\smallskip \\

Here is another way, besides non-compactness or homogeneity, to think
about the concept of ``large automorphism group''.
It turns out that one can use only the isotropy group of a single
point to characterize $B^n$ in a much more general situation. Namely,
for a complex manifold $M$ and $p\in M$, let
$I_p:=\{f\in\hbox{Aut}(M):f(p)=p\}$ be the isotropy group of $p$.

\begin{theorem}[\cite{GK3}] \sl If $M$ is a connected, non-compact
manifold of complex dimension $n$, and if there is a point $p\in M$
such that for some compact subgroup $H\subset I_p$ the set
$\{df(p):f\in H\}$ acts transitively on real tangent directions at
$p$, then $M$ is holomorphically equivalent to $B^n$.  
\end{theorem}

We would be remiss not to mention that Bland, Duchamp and Kalka \cite{BDK}
have obtained an analogue of Theorem 2.3 when the manifold
is compact.  They have also weakened the hypothesis from transitivity
on real tangent directions to transitivity on complex tangential
directions (this weakening applies both in the compact and in the
non-compact case).  Then the conclusion is that the manifold is
complex projective space.  Their techniques are different from
those in \cite{GK3}, and well worth learning. For related results see
also \cite{HO}, \cite{MN} and a discussion in \cite{GK5}.

The above results suggest several possible directions that one may
follow to endeavor to obtain characterizations for different classes of
domains with ``large''\,  automorphism group. In this survey,
we concentrate on domains whose automorphism group is non-compact, 
thus the scope of the present paper is to explore
the direction given by Theorem 2.1. 

A natural generalization of this theorem would come from replacing
the assumption of strict pseudoconvexity of $\partial D$ near a
boundary orbit accumulation point
by a weaker condition, e.g. weak
pseudoconvexity. The first results in this direction are due to 
Greene and Krantz \cite{GK4} and concern the characterization of more
general domains, namely, complex ellipsoids of the form 
$$
E_m:=\{(z_1,\dots,z_n):|z_1|^2+\dots+|z_{n-1}|^2+|z_n|^{2m}<1\}, 
$$
with $m$ a positive integer.

\begin{theorem}[\cite{GK4}] \sl  Let $D\subset {\CC}^n$ be a
bounded domain with $\hbox{Aut}(D)$ non-compact and $C^{n+1}$-smooth
boundary. Suppose that for some
boundary orbit accumulation point $p$, $\partial D$ near $p$ coincides with
$\partial E_m$ near the point $p_0=(1,0,\dots,0)\in\partial E_m$ up to a
local biholomorphism that takes $p$ into $p_0$.
Then $D$ is holomorphically equivalent to $E_m$.
\end{theorem}

If in Theorem 2.4 one allows $\partial D$ to be $C^{2n+2}$-smooth,
then the condition of local coincidence of $\partial D$ and $\partial
E_m$ up to a local biholomorphism can be replaced by the condition
that, for some local biholomorphism $f$ defined near $p$ and such that
$f(p)=p_0$, $f(\partial D)$ and $\partial E_m$ osculate to order $2m$
near $p_0$ (see \cite{GK5}). 

The proof of Theorem 2.4 uses the $C/K$-invariant that was also an
important tool in \cite{Ro}, \cite{W1}. Note that, historically, the first applications of the $C/K$-invariant to the
study of domains with non-compact automorphism group were based on
Bun Wong's results, e.g.:  
a complete hyperbolic bounded domain $D$ is biholomorphically equivalent to 
the ball if and only if there is a point $p \in D$ such 
that $[C/K](p) = 1$ \cite{W1}. A different proof of Theorem 2.4
based on the analysis of an
invariant arising from the Bergman metric can be found in
\cite{GK6}.

Kodama in \cite{Kod4} considered the more general domains
$$
E(k,\alpha):=\left\{(z_1,\dots,z_n):\sum_{j=1}^k|z_j|^2+\left(\sum_{j=k+1}^n|z_j|^2\right)^{\alpha}<1\right\},
$$
where $1\le k\le n$ and $\alpha>0$. Using methods that
avoid the $\overline{\partial}$-technique (that was needed to obtain
the main technical result of \cite{GK4}---see Lemma 4.3 there), Kodama proved a version of Theorem
2.4 for $E(k,\alpha)$ in place of $E_m$ without assuming any global
regularity of $\partial D$.  In his theorem, Kodama assumed that
$p_0\in \partial D$ and that the boundaries
$\partial D$ and $\partial E(k,\alpha)$ actually
coincide near $p_0$. In this result, one can also allow $\partial D$
and $\partial E(k,\alpha)$ to osculate near $p_0$ rather than 
literally coincide, but then a non-trivial extra condition 
on the way an orbit of $\hbox{Aut}(D)$ approaches $p_0$ is
needed.

Another generalization of Theorem 2.4 is also due to Kodama. Let
$$
E_{m_1,\dots,m_n}:=\{(z_1,\dots,z_n):
|z_1|^{2m_1}+\dots+|z_n|^{2m_n}<1\},
$$
where the $m_j$ are positive integers.

\begin{theorem}[\cite{Kod5}] \sl Let $D\subset {\CC}^n$ be a bounded
domain with non-compact automorphism group, and $p\in\partial D$ a
boundary orbit accumulation point for $\hbox{Aut}(D)$. Suppose that the
boundary of $D$ near $p$ coincides with that of $E_{m_1,\dots,m_n}$
near a point $p_0\in\partial E_{m_1,\dots,m_n}$, up to a local
biholomorphism that takes $p$ into $p_0$.  Then $D$ is holomorphically
equivalent to $E_{m_1,\dots,m_n}$.  
\end{theorem}

This result is implicit in \cite{Kod5}. There, the local equivalence
is assumed to be the identity, and the conclusion is that $D$ is
literally equal to $E_{m_1,\dots,m_n}$. But an inspection of the
proof shows that local holomorphic equivalence suffices to establish
the conclusion of global holomorphic equivalence to
$E_{m_1,\dots,m_n}$ (see \cite{GK7}). Theorem 2.5 for domains with
$C^{\infty}$-smooth boundary was obtained independently in
\cite{Ber1}. Next, in \cite{Kod7} (see also \cite{Kod6}) the version
of Theorem 2.5 as in \cite{Kod5} (stating the literal equality of the
domains) was extended to the case where the $m_i$ are arbitrary
positive real numbers (of course when the $m_i$ are not integral then the
boundary is not $C^{\infty}$-smooth). 
Further, in \cite{Kod8} this version was
proved for generalized complex ellipsoids of the form 
$$
E_{n_1,\dots,n_s;m_1,\dots,m_s}:=\left\{(z_1,\dots,z_s)\in{\CC}^{n_1}
\times\dots\times{\CC}^{n_s}:\|z_1\|^{2m_1}+\dots + \|z_s\|^{2m_s}<1
\right\}, 
$$ 
where $n_i$, $m_i$ are positive integers, $n_1+\dots+n_s=n$ and
$\|z_i\|$ denotes the ordinary norm of the vector $z_i$ in 
$\CC^{n_i}$.
In the situation where the $m_i$ are arbitrary positive real numbers,
an analogue of Theorem 2.5 for $E_{n_1,\dots,n_s;m_1,\dots,m_s}$ was
obtained in \cite{KKM}.

Further, Kim in \cite{Ki2} (see also \cite{Ki1}) obtained a result
along the lines of Theorems 2.4 and 2.5 for
domains satisfying a special local condition called {\it
Condition} $(L)$. Namely, a bounded
domain $D\subset {\CC}^n$ with non-compact automorphism group is said
to satisfy Condition $(\hbox{L})$ at a boundary orbit accumulation
point $p\in\partial D$ if $\partial D$ is real analytic near $p$, of
finite type $2k$ at $p$ in the sense of D'Angelo ($k$ is
a positive integer), and $\partial D$ near $p$ is convex up to a
local biholomorphism. The proof of Kim is in the spirit of 
the convex scaling technique due to Frankel \cite{Fr} (see also
\cite{Ki4}). As we will see later, many results below
were obtained by using a different scaling technique due to Pinchuk \cite{Pi1},
\cite{Pi3}. Since these two scaling techniques are important for
the current development of the subject, at the end of our survey we
provided a brief tutorial in the
scaling methods (for a more detailed discussion and comparison 
of these methods see \cite{Ki4}).

\section{The Bedford/Pinchuk Domains \hfill \break 
and Related Results}

In the preceding section we listed the results that extend the
Ball Characterization Theorem primarily from the point of view of the
methods and ideology suggested by its proof; in particular, we
emphasized localization principles 
that followed the work in \cite{GK4}. In this
section we turn to direct generalizations of this theorem obtained by
completely different techniques. Throughout this section all 
domains will be assumed to be smoothly bounded, i.e., bounded and
having $C^{\infty}$-smooth boundary. The first result here is due to
Bedford and Pinchuk.

\begin{theorem}[\cite{BP2}] \sl Let $D\subset{\CC}^n$ be a smoothly
bounded pseudoconvex domain of finite type with non-compact
automorphism group such that the Levi form of
$\partial D$ has no more than one zero eigenvalue at any point. Then
$D$ is holomorphically equivalent to a complex ellipsoid $E_m$ with $m$ a positive integer.  
\end{theorem}

Note that the condition on the rank of the Levi form is not a
restriction in complex dimension 2. This condition is the first step
towards allowing the domain to be weakly pseudoconvex rather than
strictly pseudoconvex; it says that
the degeneracy of the Levi form that may occur is the least possible. 
However,
in contrast with Theorem 2.1, Theorem 3.1 is essentially 
{\it non-local}.

Theorem 3.1 was first proved in \cite{BP1} for domains in ${\CC}^2$
with real analytic boundary (see also \cite{Bel3}). We note here that
real analyticity implies {\it a fortiori} the finite 
type condition  \cite{D'A2}, \cite{DF}, \cite{L}. We also mention here
that, before the paper \cite{BP2} appeared, Bell
and Catlin noticed that in \cite{BP1} real analyticity can be replaced
by the finite type condition \cite{BeCa}.

We will now give important examples of smoothly bounded domains with
non-compact automorphism group that are also due to Bedford and
Pinchuk.

\begin{example}[\cite{BP2}] \rm  Fix positive integers $m_2,\dots,m_n$
and, for a multi-index $K=(k_2,\dots,k_n)$, define its {\it weight}
by $\hbox{wt}(K)= \sum_{j=2}^n\frac{k_j}{m_j}$. Consider real
polynomials of the form 
$$
P(\tilde z,\overline{\tilde
z})=\sum_{{\rm wt}(K) = {\rm wt}(L)=1}a_{KL}{\tilde z}^K
\overline{\tilde z}^L, 
$$ 
where $\tilde z:=(z_2,\dots,z_n)$, $a_{KL}\in{\CC}$ and
$a_{KL}=\overline{a_{LK}}$. For any such polynomial we define a domain in
${\CC}^n$ by
$$
D_P:=\left\{(z_1,\tilde z):|z_1|^2+P(\tilde
z,\overline{\tilde z})<1\right\}. \eqno {(3.1)} 
$$

A domain $D_P$ of the form (3.1) is bounded if
and only if the section $D\cap\{z_1=0\}$ is a bounded domain in 
${\CC}^{n-1}$. In particular, if $D_P$ is bounded, 
then $P\ge 0$. Further, $\hbox{Aut}(D_P)$ is non-compact since it contains
the mappings
\begin{eqnarray*}
z_1 & \mapsto & \frac{z_1-a}{1-\overline{a}z_1},\\
z_j & \mapsto & \frac{(1-|a|^2)^{\frac{1}{2m_j}}z_j}
{(1-\overline{a}z_1)^{\frac{1}{m_j}}},\qquad
j=2,\dots,n,
\end{eqnarray*}
where $|a|<1$ (note that, if $D_P$ is bounded, then 
$|z_1|<1$ in $D_P$). Another way of checking the
non-compactness of $\hbox{Aut}(D_P)$ is to notice that $D_P$ is
holomorphically equivalent to the domain
$$ 
\left\{(z_1,\tilde z):\hbox{Re}z_1+P(\tilde
z,\overline{\tilde z})<1\right\}, 
$$
which is invariant under the translations
\begin{eqnarray*}
z_1 & \mapsto & z_1+it, \qquad t\in{\RR},\\
\tilde z & \mapsto &  \tilde z . \qquad \qquad \qquad \qquad \qquad 
\qquad \qquad \qquad \qquad (3.2)
\end{eqnarray*}
\null \hfill $\Box$ 
\end{example}

In their next paper \cite{BP3} Bedford and Pinchuk obtained the
following result.

\begin{theorem}[\cite{BP3}] \sl Any convex smoothly bounded domain
of finite type in ${\CC}^n$, having non-compact automorphism group, is
holomorphically equivalent to a bounded domain of the form (3.1).
\end{theorem}

The approach of Bedford and Pinchuk involves two steps. For example,
the proof of Theorem 3.1 goes as follows. In the first
step they use the method of scaling introduced
in \cite{Pi1} (see also \cite{Pi3}) to
show that the domain $D$ in consideration is holomorphically
equivalent to a domain $\Omega$ of the form
$$
\Omega=\{(z_1,\tilde z):\hbox{Re}z_1+Q(\tilde z,
\overline{\tilde z})<0\},
$$
where $Q$ is a polynomial. The domain $\Omega$ has a non-trivial
holomorphic vector field since it is invariant under translations (3.2).  
In the second step this vector field is
transported back to $D$, the result is analyzed at the parabolic fixed
point, and this information 
is used to determine the original domain. 
In the second step scaling is applied two more times. The first
scaling is needed to
show that the smallest ``weight'' involved in the vector field is
either 1 or $\frac{1}{2}$. Next, it is shown that the orbit is
well-behaved as $t\rightarrow\pm\infty$ for each of these weights, and
the final rescaling is carried out along the parabolic orbit. The case
of weight 1 is the most difficult one in the final rescaling procedure. 

There has been also certain progress, by other authors, on the first step of the above
procedure of Bedford/Pinchuk. The following completely local result
in dimension 2 (not requiring even the boundedness of the domain) was
obtained by Berteloot in \cite{Ber3} (see also \cite{BeCo},
\cite{Ber2}).

\begin{theorem}[\cite{Ber3}] \sl Let $D\subset{\CC}^2$ be a
domain, and $p\in\partial D$ a boundary orbit accumulation point for
$\hbox{Aut}(D)$. Assume that $\partial D$ is pseudoconvex and of
finite type near $p$.
Then $D$ is holomorphically equivalent to a domain of the form
$$
\{(z_1,z_2):\hbox{Re}\,z_1+P(z_2,\overline{z_2})<0\},
$$
where $P$ is a homogeneous subharmonic polynomial
without harmonic terms.
\end{theorem}

For convex domains Theorem 3.4 was recently generalized to all dimensions in
\cite{Ga}. Further, using the convex scaling technique of Frankel \cite{Fr}
 Kim in \cite{Ki2} (see also
\cite{Ki1}) obtained a related result for domains satisfying
Condition  $(\hbox{L})$ (see Section 2 for the definition).

The techniques relying on either of the two scaling principles
mentioned above (see the Appendix at the end of this paper) seem to
require the following two additional hypotheses:  pseudoconvexity (or
even convexity) and finiteness of type (or analyticity) of the
boundary. It is interesting to notice here, however, that in their
recent paper that we received while this survey was being prepared,
Bedford and Pinchuk managed to eliminate the pseudoconvexity
assumption in dimension 2 and prove the next remarkable theorem.

\begin{theorem}[\cite{BP4}] \sl Let $D\subset{\CC}^2$ be a
bounded domain with non-compact automorphism group and real analytic
boundary. Then $D$ is
holomorphically equivalent to a complex ellipsoid $E_m$
where $m$ is a positive integer.
\end{theorem}

The above results give one the hope that any smoothly bounded domain
with non-compact automorphism group should be holomorphically
equivalent to a domain of the form (3.1). Many experts believe that
this is true without extra assumptions such as the finiteness of type
and pseudoconvexity. We will now cite a result that confirms this
point of view for {\it Reinhardt domains}, i.e. domains invariant under the
rotations
$$
z_j\mapsto e^{i\phi_j}z_j,\qquad \phi_j\in{\RR},\qquad j=1,\dots,n.
$$
Note first that Reinhardt domains of the form (3.1) are given by
$$
\left\{(z_1,\tilde z): |z_1|^2+\sum_{{\rm wt}(K)=1}a_{K}{\tilde
z}^K\overline{\tilde z}^K<1\right\}, \eqno {(3.3)} 
$$
where $a_K\in{\RR}$.

\begin{theorem}[\cite{FIK2}] \sl Any smoothly bounded Reinhardt domain
in ${\CC}^n$ with non-compact automorphism group
is holomorphically equivalent to a domain of the form
(3.3), and the equivalence is
given by dilations and a permutation of coordinates.
\end{theorem}

To the best of our knowledge, Theorem 3.6 at the moment is the only
classification result for a fairly large class of domains with
non-compact automorphism group that does not require the hypotheses
of pseudoconvexity and finiteness of type. We should note, however,
that there are smoothly bounded domains with non-compact automorphism
group that are essentially non-Reinhardt. Namely, it is shown in
\cite{FIK1} that, among bounded domains of the form (3.1), there are some that are not
holomorphically equivalent to any Reinhardt domain whatsoever. The
proofs in \cite{FIK1}, \cite{FIK2} are based on the description of
the automorphism groups of bounded Reinhardt domains due
independently to Kruzhilin \cite{Kru} and Shimizu \cite{Sh}
(Kruzhilin considered the more general case of Kobayashi-hyperbolic
domains). The descriptions of Kruzhilin and Shimizu generalize that
due to Sunada for Reinhardt domains containing the origin \cite{Su}.
It is appropriate to note here that Kodama in \cite{Kod1} used the
description in \cite{Su} to prove the following

\begin{theorem}[\cite{Kod1}]  \sl Let $D\subset{\CC}^n$ be a
bounded Reinhardt domain containing the origin. Suppose that there
exists a compact subset $K\subset D$ such that $\hbox{Aut}(D)\cdot
K=D$. Then $D$ is holomorphically equivalent to a product of unit balls.
\end{theorem}

Although the situation considered in \cite{Kod1} is quite different
from that in \cite{FIK2}, the effect is essentially the same: by using the
explicit descriptions of the automorphism groups of Reinhardt domains,
one can avoid imposing any extra conditions on the boundary.
                                               
\section{The Greene/Krantz Conjecture \hfill \break
and Pseudoconvexity at Boundary \hfill \break
Orbit Accumulation Points}

As we saw in the preceding section, many of the classification
results for smoothly bounded domains with non-compact automorphism
group were proved, in particular, under the hypothesis that the
domain is of finite type. The local results in 
\cite{Ber2}, \cite{Ber3},
\cite{Ga}, \cite{Ki1}, \cite{Ki2} and local considerations in
\cite{BP1}--\cite{BP4} demonstrate the particular importance
for the boundary of the domain to be of finite type at a boundary
orbit accumulation point (note that by \cite{D'A1}---also see
\cite{D'A2}---this implies that the boundary is of finite type in
a neighborhood of the point). The Greene/Krantz conjecture states
that this geometric condition should in fact always obtain.
\medskip \\

\noindent {\bf Greene/Krantz Conjecture} {\bf([GK7])}{\sl \ \ \
Let $D\subset{\CC}^n$ 
be a smoothly bounded domain with non-compact automorphism
group. Then $\partial D$ is of finite type at any boundary orbit
accumulation point.
}
\medskip \\

The conjecture in its full generality is open. The classification in
Theorem 3.6 confirms the conjecture for Reinhardt domains, but it
would be desirable to have proofs supporting the conjecture
(even in special cases) other than those given by explicit
classification results. Below we give a theorem of such a kind due to
Kim (see also an interesting special argument presented in \cite{GK7}).

\begin{theorem}[\cite{Ki5}]  \sl Let $D\subset{\CC}^n$ be a
smoothly bounded convex domain with non-compact automorphism group, and
$p\in\partial D$. Suppose that $\partial D$ is Levi-flat 
in a neighborhood of $p$. Then $p$ is not a boundary orbit
accumulation point for $\hbox{Aut}(D)$.
\end{theorem}

Note that, as we mentioned in Section 1 above, the Levi-flatness of
$\partial D$ near $p$ is equivalent to the existence of a foliation of
$\partial D$ near $p$ by complex submanifolds of dimension
$n-1$. Thus, the relation of Theorem 4.1 to the Greene/Krantz
conjecture is that $\partial D$ does not admit such a foliation near
any boundary orbit accumulation point. However, this is, of course, a
much weaker statement compared to the conjecture itself. Also, the
conjecture is believed to be true without any extra conditions such
as the
convexity that is required in Theorem 4.1.

Another important hypothesis used in many results cited in Section 3
is the hypothesis of pseudoconvexity near a boundary orbit
accumulation point. The next theorem relates this local
pseudoconvexity to the global pseudoconvexity of the domain.

\begin{theorem}[\cite{GK6}]  \sl Let $D\subset{\CC}^n$ be a bounded domain with non-compact automorphism group, and
$p\in\partial D$ a boundary orbit accumulation point. Suppose that
$\partial D$ is $C^{\infty}$-smooth near $p$ and 
variety-free at $p$. Then local pseudoconvexity of $\partial D$ near $p$ implies that
$D$ is pseudoconvex.
\end{theorem}

Note that the variety-free assumption in Theorem 4.2 (in the case of
{\it globally} smoothly bounded domains) would follow from
the Greene/Krantz conjecture.

However, a smoothly bounded domain with non-compact automorphism
group in fact need not be globally pseudoconvex (but see
Theorem 4.4 below).  An example can be
found among the Bedford/Pinchuk domains (3.1) (see \cite{FIK2}).

\begin{example} \rm Let $D$ be the following smoothly bounded domain
(of the form (3.1)) in $\CC^3$:
$$
D:=\left\{(z_1,z_2,z_3):
|z_1|^2+|z_2|^4+|z_3|^4-\frac{3}{2}|z_2|^2|z_3|^2<1\right\}.
$$
Consider the boundary point
$p=(\frac{1}{\sqrt{2}},0,\frac{1}{\root
4\of 2})$. The complex tangent space at $p$ is
$$
\{(w_1,w_2,w_3): w_1+2^{3/4}w_3=0\} ,
$$
and the Levi form at $p$ is
$$
{\cal L}(p)(w,w) = \frac{1}{2^{\frac{3}{2}}}(-3|w_2|^2+16|w_3|^2),
$$
and thus is clearly not non-negative. Theorem 4.2 now implies that
$\partial D$ is not pseudoconvex in a neighborhood
of any boundary orbit accumulation
point (but it is {\it a fortiori} pseudoconvex {\it at} each
boundary orbit accumulation point---see Theorem 4.4 below).
\hfill $\Box$
\end{example}

We note here that in Example 4.3 the point $p$ is {\it not\/} a
boundary orbit accumulation point. In contrast, for boundary orbit
accumulation points the following general fact holds:

\begin{theorem}[\cite{GK6}] \sl Let $D\subset{\CC}^n$ be a bounded domain with non-compact automorphism group, and
$p\in\partial D$ a boundary orbit accumulation point. Suppose that
$\partial D$ is $C^2$-smooth near $p$. Then $\partial D$ is
pseudoconvex at $p$.
\end{theorem}

For domains with rough boundary a version of Theorem 4.4 remains
true.  One can say, for instance, that if $K\subset D$ is a compact
set, then its holomorphic hull $\hat K$ cannot escape to
the boundary at a boundary orbit accumulation point \cite{GK6}. For
discussions of Theorems 4.2 and 4.4 see also \cite{Kra4}.

To summarize, in this section we have discussed the two most common
hypotheses on the boundary of a smoothly bounded domain near a
boundary orbit accumulation point that occur in the literature: finiteness of type and
pseudoconvexity. The first of these is believed to always be the case
(and there are partial results to support this belief) and the second
always holds {\it at} the boundary orbit accumulation point itself;
however it is not true that the boundary
should be pseudoconvex {\it in a neighborhood} of the boundary orbit
accumulation point.  As a result of these considerations (especially
the second one), it is critical to have techniques that assume neither
pseudoconvexity nor finite type at the outset; this, in particular,
is what makes Theorem 3.5 mentioned above so important.    

\section{The Boundary Orbit Accumulation Set}

In the preceding section we dealt with individual boundary orbit
accumulation points. Here we will be interested in the collection of
{\it all} boundary orbit accumulation points, i.e. the {\it boundary orbit
accumulation set} as a whole. If $D$ is a bounded domain with
non-compact automorphism group then denote its boundary orbit
accumulation set by $S(D)$. Very little is known about the
topological and other properties of $S(D)$, except for the classes of
domains for which there exists a complete classification as in
Section 3 above. Here we give some results on $S(D)$ from \cite{IK1},
\cite{H}. The proofs use the main theorem of \cite{Bel2} and
therefore require extendability of the automorphisms to the boundary
of the domain; thus, in addition to being smoothly bounded, domains in
this section
are assumed to be pseudoconvex and of finite type (see
\cite{Kra3}).

\begin{theorem}[\cite{IK1}]  \sl  Let $D\subset{\CC}^n$ be a
smoothly bounded pseudoconvex domain of finite type with non-compact
automorphism group. Suppose that $S(D)$ contains at least three
points. Then $S(D)$ is a compact, perfect set and thus has the power of
the continuum. Moreover, in this case, $S(D)$ is either connected, or
else the number of its connected components is uncountable.
\end{theorem}

It follows from \cite{Z} that if $D$ is a 
bounded, pseudoconvex domain which is in addition {\it algebraic}, i.e.
given in the form $D=\{z\in{\CC}^n: P(z)<0\}$, where $P(z)$ is
a polynomial such that $\nabla P\ne 0$ on $\partial D$, then
the set $S(D)$ has only finitely many connected components.
Therefore, for such domains, Theorem 5.1 now implies that either $S(D)$
contains only one or two points, or $S(D)$ is connected and has
the power of the continuum.

Theorem 5.1 raises a number of natural questions: for example, can
$S(D)$ be a one- or two-point set or can $S(D)$ look like a
Cantor-type set (thus having uncountably many connected components)?
Another question is whether the set $S(D)$ is always a smooth
submanifold of $\partial D$. Note that, for instance, (the proof of)
Theorem 3.6 shows that for a smoothly bounded Reinhardt domain,
$S(D)$ is diffeomorphic to a sphere of odd dimension. The reference
\cite{GK2} gives an example of a domain with $C^{1-\epsilon}$-smooth
boundary, for which $S(D)$ has only two points.  It seems plausible
that this example can be modified, using a parabolic group of
automorphisms, so that $S(D)$ has just one point. Using similar
ideas, we also seem to be able to produce for any $k\ge 1$ a domain
$D$ with $C^k$-smooth (but not $C^{\infty}$-smooth) boundary so that
$S(D)$ has precisely two points. Indications are that the case of
finite boundary smoothness will be different from the case of
infinite boundary smoothness (see also Section 6).

If $D$ is a smoothly bounded
pseudoconvex domain of finite type, then each automorphism
of $D$ extends smoothly to the boundary.
Therefore $\hbox{Aut}(D)$ acts on the boundary, and 
the set $S(D)$ is invariant under that action.  The
following result shows that $S(D)$ is generically the
smallest invariant subset of $\partial D$.

\begin{theorem}[\cite{IK1}] \sl Let $D\subset{\CC}^n$ be a smoothly 
bounded pseudoconvex domain of finite type with non-compact
automorphism group. Suppose that $A\subset\partial D$ is non-empty,
compact and invariant under $\hbox{Aut}(D)$. Assume further that $A$
is not a one-point subset of $S(D)$. Then $S(D)\subset A$.

In particular, if $\hbox{Aut}(D)$ does not have fixed points in 
$\partial D$, then $S(D)$ is the smallest compact subset of $\partial
D$ that is invariant under $\hbox{Aut}(D)$. 
\end{theorem}

We now list several corollaries of Theorem 5.2 regarding
particular\linebreak sets $A$. Let a domain $D$
be as in the theorem. Fix $0\le k\le n-1$ and denote by
$L_k(D)$ the set of all points from $\partial D$ where the rank of
the Levi form of $\partial D$ does not exceed $k$. Clearly, each set
$L_k(D)$ is a compact subset of $\partial D$ and is invariant under
any automorphism  of $D$. Let $l_1$ denote the minimal rank of the
Levi form on $\partial D$ and $l_2$ the minimal rank of the Levi form
on $\partial D\setminus L_{l_1}(D)$.  

\begin{corollary}[\cite{IK1}, \cite{H}] \sl  Let $D$ be as in Theorem
5.2. Then either
\smallskip\\
\noindent {\bf (i)} $S(D)\subset L_{l_1}(D)$,
\smallskip \\
\noindent or
\smallskip \\
\noindent {\bf (ii)} $L_{l_1}(D)$ is a one-point subset of 
$S(D)$ and $S(D)\subset L_{l_2}(D)$.
\end{corollary}

We note that Corollary 5.3 was proved earlier by Huang \cite{H}, and its
proof in \cite{H} also relies on the paper \cite{Bel2}.

One can further endeavor to prove a property analogous
to Corollary 5.3 for the type $\tau(q)$, $q\in\partial D$, in the sense of
D'Angelo. Indeed, denote by $T_k(D)$ the set of all points
$q\in\partial D$ where $\tau(q)$ is at least $k$. We choose $t_1$ and
$t_2$ such that $T_{t_1}(D)\ne\emptyset$, $t_2<t_1$, and there exists
a point of type $t_2$ in $\partial D\setminus T_{t_1}(D)$. Since
$\tau$ is invariant under automorphisms of $D$, so is every set
$T_k(D)$. However, the sets $T_k(D)$ do not have to be closed, as the
type function $\tau$ may {\it not} be upper-semicontinuous on $\partial D$
(see e.g. an example in \cite{D'A2}, p. 136). Therefore, for the type
we only have a somewhat weaker result.

\begin{corollary}[\cite{IK1}] \sl  Let $D$ be as in Theorem 5.2. Then either
\smallskip\\ 
\noindent {\bf (i)} $S(D)\subset \overline{T_{t_1}(D)}$,
\smallskip\\
\noindent or
\smallskip\\
\noindent {\bf (ii)} $T_{t_1}(D)$ is a one-point subset of $S(D)$ and 
$S(D)\subset \overline{T_{t_2}(D)}$.
\end{corollary}

Notice that, loosely speaking, Corollaries 5.3 and 5.4 state respectively that the
rank of the
Levi form is {\it ``minimal''} and the type is {\it
``maximal''} along $S(D)$. The next corollary below states that, in this respect, 
the multiplicity function
$\mu$ (see \cite{D'A2}, p. 145 for the definition) is
analogous to the type function $\tau$. The multiplicity $\mu$ is invariant under the
extensions 
of
automorphisms to $\partial D$ and, for
$q\in\partial D$, $\tau(q)$ is finite if and only if
$\mu(q)$ is finite. In contrast with $\tau$, however, the function
$\mu$ {\it is} upper-semicontinuous on $\partial D$. 

Analogously to what we
have done above for the function $\tau$, denote by $M_k(D)$ the set of
all points $q\in\partial D$, where $\mu(q)$ is at least $k$ and choose
$m_1$ and $m_2$ such that $m_1=\hbox{max}_{q\in\partial D}\,\mu(q)$,
$m_2<m_1$, and there exists a point of multiplicity $m_2$ in $\partial
D\setminus M_{m_1}(D)$.  Because of the upper-semicontinuity and invariance
of $\mu$, each set $M_k(D)$ is a compact subset of $\partial D$
that is invariant under $\hbox{Aut}(D)$. This observation gives 
the following analogue of Corollary 5.4 for $M_{m_1}, M_{m_2}$.

\begin{corollary}[\cite{IK1}]  \sl  Let $D$ be as in Theorem 5.2. Then either
\smallskip\\
\noindent {\bf (i)} $S(D)\subset M_{m_1}(D)$,
\smallskip\\
\noindent or
\smallskip\\
\noindent {\bf (ii)} $M_{m_1}(D)$ is a one-point subset of $S(D)$ and 
$S(D)\subset M_{m_2}(D)$.
\end{corollary}

Note that one can make a statement analogous 
to Corollary 5.5 for the multitype introduced in \cite{Cat}, 
since the multitype function is upper-semicontinuous with respect to
lexicographic ordering.

It follows from Theorem 3.1 that, in complex dimension 2, for a
smoothly bounded pseudoconvex domain of finite type, the rank of the Levi form
is constant and minimal and the type is constant and maximal along
$S(D)$ (cf. Corollaries 5.3 and 5.4). Theorem 3.6 implies that this also
holds for smoothly bounded Reinhardt domains in any dimension. If we denote the minimal
rank of the Levi form by $k$, one can
see from the proof of Theorem 3.6 that for a smoothly bounded
Reinhardt domain $D$, the real dimension of any orbit of the action of
$\hbox{Aut}(D)$ on $D$ is at least $2(k+1)$. Moreover, there is
precisely one orbit of minimal dimension $2(k+1)$ (see \cite{Kra4} for a
discussion of this phenomenon). This orbit 
approaches every point of $S(D)$ non-tangentially, whereas any other
orbit approaches every point of $S(D)$ only along tangential
directions. It would be very interesting to know if similar statements hold
for more general domains. For example, the fact that there exists
an orbit that approaches $S(D)$ non-tangentially would be very
important for a proof of the Greene/Krantz conjecture (cf. 
\cite{FW1} and Theorem 4 in \cite{Ki5}). It also could
be used to show that $S(D)$ is a smooth submanifold of $\partial
D$. Generally speaking, the existence of non-tangential 
orbits to boundary orbit accumulation points---in any
(even very weak) sense---is one of the main
difficulties arising in the study of domains with non-compact
automorphism groups.

\section{More General Situations: \hfill \break
Domains with Rougher Boundary \hfill \break
and Unbounded Domains}

In this section we give certain generalizations of some of the results
mentioned above. More precisely, we will be interested in
generalizations in two directions: relaxing the condition of the
regularity of the boundary and allowing domains to be unbounded.

First we consider domains with rough boundary. Note that in Section 2 we
already mentioned the results from \cite{Kod7} and \cite{KKM} on
characterization of generalized complex ellipsoids
$E_{n_1,\dots,n_s;m_1,\dots,m_s}$ whose boundaries are not
$C^{\infty}$-smooth
if
the $m_i$ are not integers. 

The theorem below
extends Corollary 2.2 to domains with piecewise smooth boundary.
Recall that a bounded domain $D\subset\CC^n$ is said to have
{\it $C^k$-piecewise smooth} boundary, for $k\ge 1$, if $\partial D$ is a
$(2n-1)$-dimensional topological manifold and for some neighborhood
$U$ of $\partial D$ there exist real functions $\rho_j\in C^k(U)$,
$j=1,\dots,m$, such that:
\smallskip \\

\noindent {\bf (i)}  $D\cap U=\{z\in U: \rho_j(z)<0, j=1,\dots,m\}$;
\smallskip\\

\noindent {\bf (ii)}  For any subset
$\{j_1,\dots,j_r\}\subset\{1,\dots,m\}$ with $1\le j_1<\dots<j_r\le m$,
one has $d \rho_{j_1}\wedge\dots\wedge d \rho_{j_r}\ne 0$ on
$\cap_{s=1}^rS_{j_s}$, where $S_j:=\{z\in U: \rho_j=0\}$, $j=1,\dots,m$. 
\smallskip\\

The domain $D$ is said to have {\it generic $C^k$-piecewise smooth boundary}
if in the above definition one in addition has:
\smallskip\\

\noindent {\bf (iii)}  For any subset
$\{j_1,\dots,j_r\}\subset\{1,\dots,m\}$ with $1\le j_1<\dots<j_r\le m$,
one has $\partial\rho_{j_1}\wedge\dots\wedge\partial \rho_{j_r}\ne 0$
on $\cap_{s=1}^rS_{j_s}$, where $\partial$ means
differentiation only with respect to holomorphic variables.

Roughly speaking, these rather technical conditions specify
that the boundary consists of finitely many smooth pieces that have
transversal crossings.

\begin{theorem}[\cite{Pi2}] \sl  If $D\subset{\CC}^n$ is a bounded
homogeneous domain with piecewise $C^2$-smooth boundary, then $D$ is
holomorphically equivalent to a product of unit balls.
\end{theorem}

Theorem 6.1 was proved by applying the scaling method of Pinchuk that
we mentioned in Section 3 above. Note that this method also gives a
short proof of Theorem 2.1 and therefore Corollary 2.2 (see
\cite{Pi3}, \cite{Ki6} and the Appendix at the end of this survey).
Further, it was shown in \cite{CS} that for a bounded domain
with piecewise $C^2$-smooth generic boundary
such that every set $D_j:=\{z\in U: \rho_j(z)<0\}$ is strictly 
pseudoconvex, the non-compactness of $\hbox{Aut}(D)$ implies that
$D$ is in fact equivalent to $B^n$. In the case of non-tangential
approach to boundary orbit accumulation points this result was
obtained in \cite{Kod3} (see also \cite{Kod2}).

The following result is due to Kim and requires the extra hypotheses
of convexity and Levi-flatness (the latter
means that each of the sets $S_j$
from the above definition is Levi-flat).

\begin{theorem}[\cite{Ki3}] \sl Let $D\subset{\CC}^n$ be a
bounded, convex domain with piecewise $C^{\infty}$-smooth Levi-flat boundary and
non-compact automorphism group. Then $D$ is holomorphically equivalent
to the product of the unit disc and a convex domain in ${\CC}^{n-1}$.
\end{theorem}

The following local version of Theorem 6.2 is also due to Kim.

\begin{theorem}[\cite{Ki5}]  \sl Let $D\subset{\CC}^n$ be a bounded
convex domain with non-compact automorphism group. Suppose that
$\partial D$ is $C^{\infty}$-smooth and Levi-flat in a neighborhood
of some boundary orbit accumulation point. Then $D$ is
holomorphically equivalent to the product of the unit disc and a
convex domain in ${\CC}^{n-1}$.  
\end{theorem}

Note that Theorem 4.1 that we mentioned above in connection with the
Greene/Krantz conjecture is a corollary of Theorem 6.3. 

The proofs of Theorems 6.2, 6.3 rely on (an extension of) the scaling technique of
Frankel. Note that, in complex
dimension 2, these theorems give characterizations of the bidisc
$\Delta^2:=\{(z_1,z_2):|z_1|<1,|z_2|<1\}$.  Another characterization of
$\Delta^2$ is due to Bun Wong:

\begin{theorem}[\cite{W2}]  \sl  Let $D\subset{\CC}^2$ be a
bounded domain with non-compact automorphism group. Suppose that there
is a sequence $\{f_j\}\subset\hbox{Aut}(D)$ such that
\begin{eqnarray*}
&&\hbox{{\bf (i)} $W:=(\lim_{j\rightarrow\infty}f_j)(D)$ is a complex variety of
positive dimension in $\partial D$;}\\
&&\hbox{{\bf (ii)} $W$ is contained in an open subset $U\subset\partial D$
such that $U$ is} \\
&&\hbox{$C^1$-smooth and there is an 
  open set $V\subset{\CC}^2$ for
which}\\
&&\hbox{$V\cap\partial D=U$ and $V\cap D$ is convex;}\\
&&\hbox{{\bf (iii)} There exists a point $p\in D$ such that $\{f_j(p)\}$ converges to
a point} \\
&&\hbox{ $q\in W$ non-tangentially}.
\end{eqnarray*}
Then $D$ is holomorphically equivalent to $\Delta^2$.
\end{theorem}

We note here that the hypothesis {\bf (iii)} of non-tangential
convergence along some orbit is one that recurs in the literature,
but it is rather artificial.  A theorem about domains with
non-compact automorphism group should, ideally, make no hypothesis
about the way that an orbit approaches the boundary---especially a
hypothesis that is unverifiable in practice.  In fact non-tangential
approach of orbits to a boundary orbit accumulation point should be
part of the {\it conclusion} of the sorts of theorems discussed here,
not part of the hypothesis.  This hypothesis is one of the main
difficulties in problems related to domains with non-compact
automorphism groups (cf. Section 5, for example).

While this survey was being prepared we received a recent preprint of
Fu and Wong where the following result was obtained:

\begin{theorem}[\cite{FW2}] \sl Let $D\in\CC^2$ be a bounded simply-connected
domain with
generic piecewise $C^{\infty}$-smooth (but not smooth) Levi-flat boundary and
non-compact automorphism group. Then $D$ is holomorphically equivalent
to $\Delta^2$.
\end{theorem}

The next theorem deals with the case of Reinhardt domains in $\CC^2$ and is in
the spirit of Product Domain Theorems 6.2 and 6.3.

\begin{theorem}  \sl  Let $D\subset\CC^2$ be a bounded Reinhardt domain with
$C^1$-piecewise smooth (but not smooth) boundary and non-compact
automorphism group. Then $D$ is holomorphically equivalent to a
product $\Delta\times G$, where $\Delta:=\{z\in\CC:|z|<1\}$ is the
unit disc, and $G$ is either $\Delta$ or an annulus $\{1<|z|<r\}$ for
some $r>1$.
\end{theorem}

Theorem 6.6 easily follows from the proof of Theorem 3.6 in
\cite{FIK2} and holds even for domains with much rougher boundary.
The next result extends Theorem 3.6 to Reinhardt domains with
boundary of only finite smoothness.

\begin{theorem}[\cite{IK2}]  \sl   Let $D\subset{\CC}^n$ be a bounded Reinhardt
domain with $C^k$-smooth boundary, $k\ge 1$, and non-compact
automorphism group. Then, up to dilations and permutations of
coordinates, $D$ is a domain of the form
$$ 
\left\{(z_1,\dots,z_n):|z_1|^2+\psi(|z_2|,\dots,|z_n|)<1\right\},
$$ 
where $\psi(x_2,\dots,x_n)$ is a non-negative $C^k$-smooth function in 
${\RR}^{n-1}$ that is strictly positive in ${\RR}^{n-1}\setminus\{0\}$ and
such that
$\psi(|z_2|,\dots,|z_n|)$ is $C^k$-smooth
in ${\CC}^{n-1}$, and 
$$
\psi\left(t^{\frac{1}{\alpha_2}}x_2,\dots,
t^{\frac{1}{\alpha_n}}x_n\right)=t\psi(x_2,\dots,x_n) \eqno {(6.1)}
$$ 
in ${\RR}^{n-1}$ for all $t\ge 0$. Here $\alpha_j>0$, $j=2,\dots,n$, 
and each
$\alpha_j$ is either an even integer or
$\alpha_j>2k$.
\end{theorem}

In complex dimension 2, Theorem 6.7 gives the following classification:

\begin{corollary}[\cite{IK2}] \sl If $D\subset \CC^2$ is a bounded
Reinhardt domain with $C^k$-smooth boundary, $k\ge 1$, and if
$\hbox{Aut}(D)$ is non-compact, then, up to dilations and permutations
of coordinates, $D$ has the form
$$
\{|z_1|^2+|z_2|^{\alpha}<1\},
$$
where $\alpha>0$ and either is an even integer or $\alpha>2k$.
\end{corollary}

Note that, in complex dimension 3 and higher, Reinhardt domains
from Theorem 6.6 may look much more
complicated than in dimension 2 since, in contrast with the infinitely smooth case,
there is not any simple description of finitely
smooth function satisfying (6.1).

\begin{example}[\cite{IK2}] \rm  The
domain
$$
D:=\left\{|z_1|^2+|z_2|^9+|z_3|^9+\frac{1}{\log|z_3|^2-\log|z_2|^2}
\left(|z_2|^4|z_3|^5-|z_2|^5|z_3|^4\right)<1\right\}
$$
is bounded, has non-compact automorphism group (see Example 6.10) below
for a proof), and its boundary is $C^2$-smooth. The
corresponding function $\psi(|z_2|,|z_3|)$ possesses weighted
homogeneity property (6.1) with $\alpha_2=\alpha_3=9$. \\
\null \hfill $\Box$
\end{example}

Along the lines of Theorem 6.7, one can consider the following examples of
domains with non-compact automorphism group and
$C^k$-smooth boundary, $k\ge 1$, that are not necessarily
Reinhardt.

\begin{example}[\cite{IK2}] \rm Consider the domain
$$ 
\left\{(z_1,\dots,z_n):|z_1|^2+\psi(z_2,\dots,z_n)<1\right\},\eqno {(6.2)}  
$$ 
where $\psi(z_2,\dots,z_n)$ is a $C^k$-smooth function on
${\CC}^{n-1}$ and 
$$
\psi\left(t^{\frac{1}{\alpha_2}}z_2,\dots,
t^{\frac{1}{\alpha_n}}z_n\right)=|t|\psi(z_2,\dots,z_n)\eqno {(6.3)}
$$  
in ${\CC}^{n-1}$ for all $t\in{\CC}\setminus\{z:\Re z<0\}$. Here
$\alpha_j>0$, $j=2,\dots,n$, and
$t^{\frac{1}{\alpha_j}}=e^{\frac{1}{\alpha_j}(\log|t|+i{\rm arg}\,t)}$,
for $t\ne 0$ and $-\pi<\hbox{arg}\,t<\pi$. Also, to guarantee that
the domain given in (6.2) is bounded, one can assume that
the domain in $\CC^{n-1}$   
$$
\left\{(z_2,\dots,z_n):
\psi(z_2,\dots,z_n)<1\right\}  
$$ 
is bounded. 

For any domain $D$ of the form (6.2), $\hbox{Aut}(D)$ is indeed 
non-compact, since it contains the subgroup
\begin{eqnarray*}
z_1& \mapsto & \frac{z_1-a}{1-\overline{a}z_1},\\
z_j& \mapsto & \frac{(1-|a|^2)^{\frac{1}{\alpha_j}}z_j}
{(1-\overline{a}z_1)^{\frac{2}{\alpha_j}}},\qquad j=2,\dots,n,
\end{eqnarray*}
where $|a|<1$.
\hfill $\Box$
\end{example}

If $n=2$ then, by differentiating both parts of (6.3) with respect to $t$ 
and $\overline{t}$ and setting $t=1$, we obtain that
$\psi(z_2)=c|z_2|^{\alpha}$, with $c>0$. Therefore, for $n=2$, the domain
(6.2) is equivalent to a domain of the form 
$$
\{(z_1,z_2):|z_1|^2+|z_2|^{\alpha}<1\}
$$
which is Reinhardt. However, as
examples in \cite{FIK1} show, there exists a bounded domain in
${\CC}^2$ with non-compact automorphism group whose boundary 
is {\bf (i)}  real analytic at all points except
one, {\bf (ii)} $C^{1,\beta}$-smooth at the exceptional point for some
$0<\beta<1$, and that is biholomorphically inequivalent
to any Reinhardt domain and thus to any domain of the form (6.2).
It would be interesting to construct such examples for the case of
$C^k$-smooth boundaries, $k\ge 2$ and in dimensions $n\ge 2$. The domains (6.2) seem to 
be reasonable generalizations of the
Bedford/Pinchuk domains (3.1) to the case of finitely smooth boundaries.

We turn now to the case of unbounded domains. Note that in
Section 3 we already gave some classification results that hold for
unbounded domains because of their completely local nature (see e.g. Theorem
3.4 and \cite{Ga}).  Another local result that we
mention here is due to Efimov \cite{E} and generalizes Theorem 2.1 to the case
of unbounded domains (note that Theorem 2.1, being local, still
requires the domain to be bounded).

\begin{theorem}[\cite{E}]  \sl  
Let $D\subset{\CC}^n$ be a domain
(not necessarily bounded), and $p\in\partial D$ a boundary orbit
accumulation point for $\hbox{Aut}(D)$. Assume that $\partial D$ is
$C^2$-smooth and strictly pseudoconvex near $p$. Then $D$ is
holomorphically equivalent to $B^n$. 
\end{theorem}

The next theorem is not local, and the domain is assumed to be
Kobayashi-hyperbolic.

\begin{theorem}[\cite{IK3}]  \sl  Let $D\subset{\CC}^2$ be
a hyperbolic Reinhardt domain with $C^k$-smooth
boundary, $k\ge 1$, and let $D$ intersect at least one of
the coordinate complex lines $\{z_j=0\}$, $j=1,2$. Assume
also that $\hbox{Aut}(D)$ is non-compact. Then $D$ is
holomorphically equivalent to one of the following
domains:
\begin{eqnarray*}
&&\hbox{\bf (i)}\,\left\{(z_1,z_2):|z_1|^2+|z_2|^{\alpha}<1\right\},\\
&&\qquad \hbox{where either $\alpha<0$, or
$\alpha=2m$ for some $m\in{\NN}$, or $\alpha>2k$;}\\
&& \hbox{ \ } \\
&&\hbox{\bf (ii)}\,\left\{(z_1,z_2):|z_1|<1,(1-|z_1|^2)^{\alpha}<|z_2|<
R(1-|z_1|^2)^{\alpha}\right\},\\
&& \qquad \hbox{where $1<R\le\infty$ and $\alpha<0$;}\\
&& \hbox{ \ } \\
&&\hbox{\bf (iii)}\,\left\{(z_1,z_2):e^{\beta|z_1|^2}<|z_2|<Re^{\beta|z_1|^2}\right\},\\
&& \qquad \hbox{where $1<R\le\infty$, $\beta\in{\RR}$, $\beta\ne
0$, and, if $R=\infty$, $\beta>0$.} 
\end{eqnarray*}

If $k<\infty$ and $\partial D$ is not
$C^{\infty}$-smooth, then $D$ is holomorphically
equivalent to a domain of the form {\bf (i)} with
$\alpha\ne 2m$ for any
$m\in{\NN}$ and $\alpha>2k$.

In case {\bf (i)} the equivalence is given by dilations and
a permutation of the coordinates; in cases {\bf (ii)} and
{\bf (iii)} the equivalence is given by a mapping of the form
\begin{eqnarray*}
z_1&\mapsto& \lambda z_{\sigma(1)}
z_{\sigma(2)}^a,\\ 
z_2&\mapsto& \mu z_{\sigma(2)}^{\pm 1},
\end{eqnarray*}
where $\lambda,\mu \in {\CC}^*$, $a\in{\ZZ}$ and
$\sigma$ is a permutation of $\{1,2\}$.
\end{theorem}

Note that in Theorem 6.12 we do not assume the existence of a finite
boundary orbit accumulation point (of course the domain may be
unbounded)
which is an important hypothesis in \cite{Ber3}, \cite{Ga}. The condition for the domain to
intersect a coordinate complex line gives that $\hbox{Aut}(D)$ has
only finitely many connected components and therefore the non-compactness
of $\hbox{Aut}(D)$ is equivalent to the non-compactness 
of its identity component, to which the description in \cite{Kru} applies.
It seems that without this condition there is not any
reasonable classification, since one can produce many ``exotic''
hyperbolic domains for which the identity component of the
automorphism group is compact whereas the whole group is non-compact and
has infinitely many connected components. Domains with such a structure
of the automorphism groups seem to be intractable. To support this
claim we give one ``exotic'' example below.

\begin{example}[\cite{IK3}] \rm Consider the Reinhardt domain $D\subset
{\CC}^2$ given by
$$
D:=\left\{(z_1,z_2):\sin\left(\log\frac{|z_1|}{|z_2|}\right)<
\log |z_1
z_2|<\sin\left(\log\frac{|z_1|}{|z_2|}\right)+\frac{1}{2}\right\}.
$$
The boundary of $D$ is clearly $C^{\infty}$-smooth. The
group $\hbox{Aut}(D)$ is not compact since it contains
all the mappings 
\begin{eqnarray*}
z_1& \mapsto & e^{\pi k}z_1,\\
z_2& \mapsto & e^{-\pi k}z_2,
\end{eqnarray*}
for $k\in{\ZZ}$.

To see that $D$ is hyperbolic, consider the mapping
$f: D\rightarrow {\CC}$, $f(z_1,z_2)=z_1 z_2$. It
is easy to see that $f$ maps $D$ onto the annulus
$A:=\left\{z \in \CC:e^{-1}<|z|<e^{\frac{3}{2}}\right\}$, which is a
hyperbolic domain in ${\CC}$. The annuli
\begin{eqnarray*}
A_1&:=&\left\{z \in \CC:e^{-\frac{1}{4}}<|z|<e^{\frac{1}{2}}\right\}, \\ 
A_2&:=&\left\{z \in \CC:e^{-1}<|z|<e^{-\frac{1}{8}}\right\}, \\
A_3&:=&\left\{z \in \CC:e^{\frac{1}{4}}<|z|<e^{\frac{3}{2}}\right\}
\end{eqnarray*}
obviously cover $A$, and each of the inverse images
$D_j=f^{-1}(A_j)$, $j=1,2,3$, is hyperbolic since
$D_j$ is contained in a union of bounded
pairwise non-intersecting domains. It then follows (see
\cite{PoSh}) that $D$ is hyperbolic.
\hfill $\Box$
\end{example}

The following example
suggests that, in complex dimension $n\ge 3$, an explicit
classification result in the hyperbolic case---analogous to 
Theorem 6.12---does not exist; in fact it does
not exist if we do not impose extra conditions on the
domain, even if the domain contains the origin. We note here that, to
obtain the finiteness of the number of connected components of
$\hbox{Aut}(D)$ for a hyperbolic Reinhardt domain $D\subset\CC^n$, one
needs the assumption that $D$ intersects at least $n-1$ coordinate
hyperplanes \cite{IK3}, which certainly holds for domains containing
the origin. Therefore
the problem suggested by the following example
is of a different kind compared to the one arising from Example 6.12
and is specific for dimensions $n\ge 3$. 
\medskip

\begin{example}[\cite{IK3}] \rm Consider the domain $D\subset{\CC}^3$ given by
\begin{eqnarray*}
D&:=& \Biggl\{(z_1,z_2,z_3): \phi(z):=  \\
 && \qquad |z_1|^2 + (1-|z_1|^2)^2|z_2|^2\rho\Bigl(|z_2|^2(1-|z_1|^2) , 
         |z_3|^2(1-|z_1|^2)\Bigr) \\ 
&&  \qquad \qquad
 +(1-|z_1|^2)^2|z_3|^2 - 1 < 0 \Biggr\}, \qquad \qquad 
\qquad \qquad \qquad {(6.4)}
\end{eqnarray*}
where $\rho(x_1,x_2)$ is a
$C^{\infty}$-smooth function on ${\RR}^2$ such that
$\rho(x_1,x_2)>c>0$ everywhere, and the partial
derivatives of $\rho$ are  non-negative for $x_1,x_2\ge 0$.

To show that $\partial D$ is smooth, we calculate
\vfill
\eject

\begin{eqnarray*}
\frac{\partial\phi}{\partial
z_1}&=&\overline{z_1}\Biggl(1-(1-|z_1|^2)\Bigl(2|z_2|^2\rho+
(1-|z_1|^2)|z_2|^4\frac{\partial\rho}{\partial
x_1}\\
&& \qquad \qquad \qquad
+(1-|z_1|^2)|z_2|^2|z_3|^2\frac{\partial\rho}{\partial
x_2}+2|z_3|^2\Bigr)\Biggr),\\
\frac{\partial\phi}{\partial
z_2}&=&(1-|z_1|^2)^2\overline{z_2}\left(\rho+(1-|z_1|^2)|z_2|^2
\frac{\partial \rho}{\partial x_1}\right),\\
\frac{\partial\phi}{\partial
z_3}&=&(1-|z_1|^2)^2\overline{z_3}\left((1-|z_1|^2)|z_2|^2
\frac{\partial\rho}{\partial x_2}+1\right). \qquad \qquad {(6.5)}
\end{eqnarray*}
It follows from (6.5) that not all the partial derivatives of
$\phi$ can vanish simultaneously at a point of
$\partial D$. Indeed, if $\frac{\partial\phi}{\partial
z_3}(p)=0$ at some point $p\in\partial D$ then, at $p$,
either $|z_1|=1$ or $z_3=0$. If $|z_1|=1$, then clearly
$\frac{\partial\phi}{\partial z_1}(p)\ne 0$. If $|z_1|\ne
1$, $z_3=0$, and, in addition,
$\frac{\partial\phi}{\partial z_2}(p)=0$, then $z_2=0$,
and therefore $|z_1|=1$, which is a contradiction.
Therefore, $\partial D$ is $C^{\infty}$-smooth.

To show that $D$ is hyperbolic, consider the holomorphic
mapping defined by $f(z_1,z_2,z_3)=z_1$ from $D$ into ${\CC}$.
Clearly $f$ maps $D$ onto the unit disc
$\Delta$, which is a hyperbolic domain in
${\CC}$. Further, the discs $\Delta_r:=\{z:|z|<r\}$ for
$r<1$ form a cover of $\Delta$, and $f^{-1}(\Delta_r)$
is a bounded open subset of $D$ for any such $r$. Thus,
as in Example 6.13 above, we see that $D$ is
hyperbolic (see \cite{PoSh}).

Further, $\hbox{Aut}(D)$ is non-compact since it contains
the automorphisms
\begin{eqnarray*}
z_1 & \mapsto & \frac{z_1-a}{1-\overline{a}z_1},\\
z_2 & \mapsto & \frac{(1-\overline{a}z_1)z_2}{\sqrt{1-|a|^2}},\qquad
\qquad \qquad \qquad \qquad \qquad {(6.6)}\\
z_3 & \mapsto & \frac{(1-\overline{a}z_1)z_3}{\sqrt{1-|a|^2}},
\end{eqnarray*}
for $|a|<1$.
\hfill $\Box$
\end{example}

Examples similar to Example 6.14 can be constructed in any complex
dimension $n\ge 3$. They indicate that, most probably,
there is no reasonable classification of smooth hyperbolic Reinhardt
domains with non-compact automorphism group for $n\ge 3$ even in the
case when the domains contain the origin. Indeed, in Example 6.14 we
have substantial freedom in choosing the function $\rho$. 
We note
that the boundary of domain (6.4) contains the complex hyperplane
$z_1=\alpha$ for any $|\alpha|=1$. It may happen that, by imposing
the extra condition of the finiteness of type on the boundary of the
domain, one would eliminate the difficulty arising in Example 6.14 and
obtain an explicit classification. It also should be observed that
any point of the boundary of domain (6.4) with $|z_1|=1$, $z_2=z_3=0$
is a boundary orbit accumulation point for $\hbox{Aut}(D)$ (see (6.6));
therefore, it is plausible that one needs the finite type condition
only at such points (cf. the Greene/Krantz conjecture for the bounded case).

\section{Concluding Remarks}

The study of automorphism groups has considerable intrinsic interest,
and also has roots in several of the major themes of twentieth
century mathematics.  Because domains in higher dimensions are
generically biholomorphically distinct, it is natural to seek
some unifying properties that they enjoy.  The automorphism
group provides one such natural set of ideas.

The program we have described suggests that considerable progress
has been made in understanding domains
of finite type with ``large'' automorphism group. The Greene/Krantz conjecture, which at this
point in time appears likely to be true, suggests that finite
type domains are the only ones that require study.

However, it should be borne in mind that these last remarks apply
only to smoothly bounded domains.  Evidence suggests that each
boundary smoothness class $C^k$ has different automorphism group
phenomena, and that the picture becomes more and more complicated
as $k$ becomes smaller.  In particular, for domains with fractal
boundary almost nothing is known (and the self-similarity of a fractal
boundary suggests that this case is of particular interest for
automorphism group symmetry).  We look forward to new insights in
the future, some perhaps inspired by the present article.
\bigskip
\bigskip\\

\begin{center}
----------------------------------------------------------------------------------------------- \\
----------------------------------------------------------------------------------------------- \\
\end{center}
\vspace*{.13in}

\newpage

\begin{center}
{\Large \bf Appendix on the Scaling Methods}
\end{center}

We now sketch the key ideas in the methods of scaling and some of their
applications to the study of domains with non-compact automorphism
groups. We begin
with the method originated by S. Pinchuk in the late 1970's \cite{Pi1}.
This discussion will lead to a proof of the Ball Characterization 
Theorem (Theorem 2.1). We then conclude the Appendix with
an outline of Frankel's scaling technique. Our discussion of scaling mainly follows
the exposition in \cite{Ki4}.

In the discussion of Pinchuk's method, for simplicity, we restrict attention to
scaling of strongly pseudoconvex domains. This will convey
the main ideas without the added baggage that treating
finite type points would entail. It should be clearly understood,
however, that scaling is of greatest importance in the weakly
pseudoconvex case because it is virtually the only technique
available in that setting.

Fix a smoothly bounded domain $D$ with strongly pseudoconvex
boundary point $q$.  We assume that $q$ is a boundary orbit accumulation
point for the action of the automorphism group on $D$.  Therefore
there are a point $p\in D$ and a sequence of automorphisms $\{f_j\}\subset\Aut(D)$
such that $f_j(p) \rightarrow q$ as $j\ra \infty$. We may apply a quadratic
holomorphic polynomial change of coordinates so that $q$ is
mapped into the origin and there is a ball $U$ centered at 
the origin such that 
$U \cap \partial D$ 
is strongly {\it convex} (see Narasimhan's Lemma
in \cite{Kra3}). Denote
$\tilde z:=(z_2,\dots,z_n)$, so that $z:=(z_1,\dots,z_n)=(z_1,\tilde z)$.
Now a simple holomorphic change of coordinates (we denote it by $F$)
allows us to write a defining function on the set $U \cap D$
(with a possibly smaller ball $U$) as
$$
\rho(z_1,\tilde z):=\Re z_1 + ||\tilde z||^2 + o\bigl( |\Im z_1| +
||\tilde z||^2\bigr).
\eqno (A.1) 
$$
It follows then that $\partial D$ is variety-free at $q$. Now a
simple normal family argument implies the following result (see \cite{Kra3} for details):
\smallskip \\

\noindent {\bf Lemma A.1} {\sl Let notation be as above.
Then there is a subsequence of $\{f_j\}$ that converges to the
constant mapping $q$ uniformly on compact subsets of $D$.}
\vspace*{.1in}

Define $p^j = f_j(p)$ for each $j$. Of course $p^j \ra p$
as $j \ra \infty$. Set $p^j = (p^j_1, \dots, p^j_n)$.
For each $j$, we construct a holomorphic change of variables
as follows:
$$
\left \{ \begin{array}{lcl}
  \zhat_1 & = & e^{i\theta_j}z_1 - p_j^* - \sum_{m=2}^n a_m (z_m - p^j_m), \\
  \tilde{\zhat} & = & \tilde{z} - \tilde{p^j}.
\end{array}
\right.
\eqno (A.2)
$$
Here $\theta_j\in\RR$ and $p_j^*, a_m\in\CC$ are selected so that in
the coordinates $\zhat:=(\zhat_1,\dots,\zhat_n)$ one has:
\begin{itemize}
\item $(0,\dots,0)\in\partial D$;
\item $p^j=(-\e_j , 0, \dots,0)$, $\e_j > 0$, for each $j$;
\item The tangent plane to $\partial D$ at $(0,\dots,0)$ is given
by $\{z: \Re z_1 = 0\}$.
\end{itemize}

In the $\zhat$-coordinates the defining function in equation $(A.1)$
is given by
$$
\widehat{\rho}_j(\zhat):=\widehat{c}_j \Re \left (
    \zhat_1 + \sum_{m=1}^n A_m^j \zhat_m^2 \right ) +
              \sum_{k,m = 1}^n B_{km}^j \zhat_k \overline{\zhat_m} + 
E_j(\zhat),
$$
where $E_j(\zhat)=o(|\Im \zhat_1| + ||\tilde{\zhat}||^2)$ and the coefficients
of the quadratic terms converge to the coefficients of the 
corresponding quadratic terms for the defining function $\rho$ 
in $(A.1)$.  Furthermore, $\widehat{c}_j \ra 1$ as $j \ra \infty$.

Now we come to the heart of the scaling process.  Thus far we have
been normalizing coordinates so that the scaling can be performed
in a natural manner.  The motivation for the scaling that we do
is as follows:  the natural geometry of a strongly pseudoconvex
point is parabolic in nature.  This can be seen by examining
the boundary behavior of the Carath\'{e}odory or Kobayashi metrics
(see \cite{Kra3}), but can be also seen in a more elementary fashion
by examining, for instance, the defining function in $(A.1)$.  
We see that an arbitrary strongly pseudoconvex point can be
viewed as a perturbation of the domain 
$$
\widehat{D}:=\{z\in\CC^n:
\widehat{\rho}:=\Re z_1 + ||\tilde z||^2 < 0\}.\eqno (A.3)
$$

A moment's thought reveals this last
domain to be holomorphically equivalent to the unit ball $B^n$
(see e.g. \cite{Ru}). And the
parabolic
nature of the boundary is self-evident from comparing the roles of
$\Re z_1$ and $\tilde z$ in
the defining function $\widehat{\rho}$. 

Having said all this, we now set
$$
\begin{array}{lcl}
   \displaystyle    z'_1 &=& \displaystyle \frac{\widehat{z}_1}{\epsilon_j},\\
                   \hbox{ \ \ } & \hbox{ \ } & \hbox{ \ } \\
   \displaystyle   \tilde{z'} &=& \displaystyle \frac{\tilde{\widehat{z}}}{\sqrt{\epsilon_j}}, \\
\end{array} \eqno (A.4)
$$
with $\e_j$ defined by $(A.2)$. Given that a strongly pseudoconvex point is
nearly like a ball, what we are doing is scaling that ball up
to have radius about 1.  But the magnitude of the scaling depends
on the normal distance of $p_j$ to the boundary.   

Let $D_j$ denote the image of
$D\cap U$ under the composition of $F$, mapping $(A.2)$ and mapping $(A.4)$.
Taking into account the fact that $\e_j \ra 0$ as $j \ra \infty$,
we may write (dropping primes) the defining function for $D_j$ as
$$
\rho_j(z)=c_j \Re \left ( z_1 + \sum_m A_m^j z_m^2 \right)
            + \sum_{k,m} B_{km}^j z_k \zbar_m
            + \e_j^{-1} E_j(\e_j z_1, \sqrt{\e_j} \tilde z) .
$$

As $j \ra \infty$, we see that the ``limiting defining function''
is then the function $\widehat{\rho}$\/ from $(A.3)$.
Then the domains $D_j$ converge (in the sense of Hausdorff set
convergence) to the limiting domain $\widehat{D}$.

Now the crux of the matter is this: combining our various coordinate
changes, we see that we have constructed, for each $j$, 
a biholomorphic mapping
$$
g_j: U \cap D \longrightarrow D_j.
$$
For any 
compact subset $K\subset\widehat{D}$, we have $K\subset D_j$ for $j$
sufficiently large and thus $G_j:=f_j^{-1} \circ g_j^{-1}$ is defined
on $K$ for large $j$. Since $D$ is bounded and $K$ was an arbitrary
compact  subset of $\widehat{D}$, a subsequential limit yields a holomorphic
mapping
$$
g: \widehat{D} \longrightarrow {\overline D}.
$$
On the other hand, by Lemma A.1, passing if necessary to a subsequence, we also know that any
compact subset of $D$ is mapped to $D\cap U$ under $f_j$ for large
$j$. Therefore, $G_j^{-1}=g_j\circ f_j$ is defined on any compact subset of $D$
for $j$ large enough. 

Using these two facts, it is possible to prove that the limit mapping
$g$ is in fact a
biholomorphism from $\widehat{D}$ onto $D$ (see e.g. \cite{Ki4}). Since $\widehat{D}$ is holomorphically equivalent to
$B^n$, so is $D$.  This concludes
the proof of the Ball
Characterization Theorem by scaling. 

At the level of strongly pseudoconvex domains, the scaling technique
is largely formalistic. In the case of weakly pseudoconvex domains of
finite type the argument just presented is only the beginning of the
proof. The difficulty in this case is that the limit domain $\widehat
D$ is not so easily found as for strictly pseudoconvex domains. To
determine $\widehat{D}$ one needs
an argument that involves further applications
of the scaling process \cite{BP1}--\cite{BP4}.

With this last thought in mind, we now say just a few words
about a scaling technique introduced by Frankel \cite{Fr}. It has proved
to be important because, in the case when the domains under consideration
are {\it convex}, the delicate limiting arguments described above are easier to handle. Note that in the proof of
Theorem 3.3 convexity also helps to make scaling 
arguments---based on Pinchuk's method--- easier (see \cite{BP3}).

Now let $D\subset\CC^n$ be a bounded, convex domain.  Suppose, as before, that
there are a point $p \in D$ and a sequence of automorphisms $\{f_j\}$ of $D$
such that $f_j(p) \ra q \in \partial D$.  Consider the mappings
$$
\omega_j(z) = \bigl [ \partial f_j(p)\bigr ]^{-1} (f_j(z)-f_j(p)) .
$$
where $\partial f_j$ is the holomorphic Jacobian matrix of $f_j$. The
central point of the scaling procedure is the following result of Frankel.
\medskip\\

\noindent {\bf THEOREM A.2 (\cite{Fr})} {\sl Let notation be as
above. Then

\noindent {\bf (i)} $\{\omega_j\}$
is a normal family (i.e. every subsequence of $\{\omega_j\}$ has a
subsequence that uniformly converges on compact subsets of $D$);

\noindent {\bf (ii)} Every subsequential limit of $\{\omega_j\}$ is a
holomorphic embedding of $D$ into $\CC^n$.}
\medskip\\

The following version of the above result of Frankel is due to
Kim.
\medskip\\

\noindent {\bf PROPOSITION A.3 (\cite{Ki4})} {\sl Let notation be as
above. Suppose that $\partial D$ is variety-free at $q$. Then, by
passing to a subsequence of $\{f_j\}$ if necessary, one can construct a sequence $\{q_j\}\subset\partial D$, $q_j\ra q$ as
$j\ra\infty$, such that
\smallskip\\

\noindent {\bf (i)} The mappings
$$
\sigma_j(z) = \bigl [ \partial f_j\bigr(p) ]^{-1} (f_j(z)-q_j)
$$
form a normal family;
\smallskip\\

\noindent {\bf (ii)} Every subsequential limit of $\{\sigma_j\}$ is a
holomorphic embedding of $D$ into $\CC^n$.}
\medskip\\

Note that
Theorem A.2 and Proposition A.3 do not require any
regularity of $\partial D$. The sequence of scaled domains that one
has to consider is then the sequence
$\{\sigma_j(D)\}$. Further, as in Pinchuk's method above, one has to understand what the limit
domain is and why it is holomorphically equivalent to $D$, and this is
where the regularity of $\partial D$ becomes important.

\bigskip

{\obeylines
Centre for Mathematics and Its Applications 
The Australian National University 
Canberra, ACT 0200
AUSTRALIA 
E-mail address: Alexander.Isaev@anu.edu.au
\smallskip
Department of Mathematics
Washington University, St.Louis, MO 63130
USA 
E-mail address: sk@math.wustl.edu}

\end{document}